\documentclass{amsart}
\usepackage{amsfonts}
\usepackage{fancyvrb}
\newtheorem{thm}{Theorem}

\newtheorem{rem}{Remark}
\newtheorem{exa}{Example}

\usepackage{color}
\definecolor{brown}{rgb}{0.8,0.6,0.3}
\definecolor{dgreen}{rgb}{0.2,0.4,0.3}

\usepackage[pdfauthor={Richard J. Mathar}, colorlinks=true,urlcolor=blue,linkcolor=red,citecolor=brown,
pdfkeywords={{generating function},{P-finite},{D-finite}},bookmarksopen=true,pdfpagelayout=OneColumn]{hyperref}

\begin{document}

\title[D-finite Generating Functions]{P-finite Recurrences From Generating Functions with Roots of Polynomials}

\author{Richard J. Mathar}
\urladdr{https://www.mpia-hd.mpg.de/~mathar}
\address{Max-Planck Institute of Astronomy, K\"onigstuhl 17, 69117 Heidelberg, Germany}
\email{mathar@mpia-hd.mpg.de}


\date{\today}

\begin{abstract}
We derive the P-finite recurrences for classes of sequences with 
ordinary generating function containing roots of polynomials. The focus is on
establishing the D-finite differential equations such that the
familiar steps of reducing their power series expansions apply.
\end{abstract}

\maketitle
\section{Aim}
The manuscript derives P-finite recurrences for some classes of generating functions $g(x)$
which have a self-replicating property under differentiation, mainly involving roots
and exponentials of polynomials of $x$.
The standard application of the recurrences are
\begin{itemize}
\item
numerical generation of a deep list of the expansion coefficients
with constant time and with constant requirements on memory for each additional term.
[George Fischer's work on efficiently implementing a large set of holonomic sequences
of the Online Encyclopedia of Integer Sequences (OEIS) was in fact the main inspiration
of this work, see \url{https://github.com/archmageirvine/joeis}.]
\item
asymptotic estimates for large $n$ \cite{WimpJMAA111}\cite[Thm. VII.10]{Flajolet},
where the generating functions cover statistics of walks, biological cycles and similar
growths to higher generations.
\item
reverse engineering of generating functions if the P-finite recurrences 
match one of the emerging types.
\end{itemize}

\section{Inverse Root}\label{sec.iroot} 
This section is a tutorial which demonstrates
the mechanism of coefficient shifts for differentiation of power series.  
Let $g(x)$ be the (ordinary) generating function of a sequence which is an inverse $r$-th
root of a polynomial $p(x)$,
\begin{equation}
g(x) \equiv \sum_{n\ge 0} g_n x^n = \frac{1}{\sqrt[r]{p(x)}}, \quad r\neq 0.
\label{eq.gdef}
\end{equation}
\begin{equation}
p(x) \equiv \sum_{n=0}^{\deg p} p_n x^n.
\end{equation}
Then the chain rule of differentiation yields a first derivative of  \eqref{eq.gdef},
\begin{equation}
g'(x) = -\frac{1}{r}\frac{p'(x)}{p^{1+1/r}(x)}
 = -\frac{1}{r}\frac{p'(x)}{p(x)}g(x)
.
\end{equation}
Obviously $g(x)$ is a D-finite function \cite{StanleyEJC1,LipshitzJA122}:
\begin{equation}
r p(x) g'(x) + p'(x)g(x) =0.
\end{equation}
The P-finite recurrence is derived by insertion of the two power series:
\begin{equation}
r\sum_{i=0}^{\deg p} p_i x^i \sum_{j\ge 1} j g_j x^{j-1} +
\sum_{i=1}^{\deg p} i p_i x^{i-1} \sum_{j\ge 0} g_j x^j =0.
\end{equation}
Resummation with $k\equiv i+j-1$ yields
\begin{multline}
r\sum_{k=0}^{\deg p -1} \sum_{j=1}^{k+1} p_{k+1-j} j g_j x^k
+r\sum_{k\ge \deg p} \sum_{j=k+1-\deg p}^{k+1} p_{k+1-j} j g_j x^k \\
+ \sum_{k=0}^{\deg p-2} \sum_{j=0}^k (k+1-j) p_{k+1-j} g_j x^k 
+ \sum_{k\ge \deg p-1} \sum_{j=k+1-\deg p}^k (k+1-j) p_{k+1-j} g_j x^k =0.
\end{multline}
Comparison of coefficients for $x^k$, $k\ge \deg p$, on both sides gives
\begin{equation}
r \sum_{j=k+1-\deg p}^{k+1} p_{k+1-j} j g_j 
+ \sum_{j=k+1-\deg p}^k (k+1-j) p_{k+1-j} g_j =0,\quad k\ge \deg p.
\end{equation}
Setting $n\equiv k+1$ gives
\begin{equation}
r \sum_{j=n-\deg p}^{n} j p_{n-j} g_j 
+ \sum_{j=n-\deg p}^{n-1} (n-j) p_{n-j} g_j =0,\quad n\ge \deg p+1;
\end{equation}
\begin{equation}
r n p_0 g_n 
+ \sum_{j=n-\deg p}^{n-1} [rj +n-j]p_{n-j} g_j 
=0,\quad n\ge \deg p+1.
\end{equation}
Setting $j=n-t$ gives
\begin{equation}
r n p_0 g_n 
+\sum_{t=1}^{\deg p} (rn-rt+t)p_t g_{n-t} 
=0,\quad n\ge \deg p+1.
\end{equation}
As an extension of Noe's \cite{NoeJIS9}
recurrences we find:
\begin{thm}
The coefficients of the generating function \eqref{eq.gdef} obey the P-finite recurrence
\begin{equation}
\sum_{t=0}^{\deg p} (rn-rt+t)p_t g_{n-t} 
=0,\quad n\ge {\deg p}+1.
\end{equation}
\end{thm}

\begin{rem}
For $r=1$ this equation may
be divided through a common factor $n$, 
\begin{equation}
\sum_{t=0}^{\deg p} p_t g_{n-t} =0,\quad n\ge {\deg p}+1, \quad r=1,
\end{equation}
which is a recurrence with constant coefficients $p_t$.
We recover the well known result:
The sequences with a rational generating function $p(x)/q(x)$
have recurrences where the D-finite equation does not contain derivatives of $g(x)$,
so the $g(n)$ obey C-finite recurrences.
\end{rem}

\begin{rem}
Our formulas for recurrences normalize the representation
by using indices $g_{n-t}$, $t\ge 0$, because that is the most
readable convention while implementing computer programs that
derive sequences from lower-order terms \cite{NemesJSC20,ReutenauerEJC19}.
\end{rem}

\begin{exa}
Examples of \eqref{eq.gdef}, square roots $r=2$, degree $\deg p=1,2$, including
the OEIS numbers of the coefficient sequences \cite{sloane}

\begin{tabular}{ll}
$p(x)$ & \\
\hline
$1-2x-3x^2$ & A002426 \\
$1-6x-3x^2$ & A122868 \\
$1-6x+x^2$ & A001850 \\
$1-6x+5x^2$ & A026375 \\
$1-4x-4x^2$ & A006139 \\
$1-4x$ & A000984 \\
$1-2x+5x^2$ & A098331 \\
$1-4x^2$ & A126869 \\
\end{tabular}

\end{exa}

\begin{exa}
Examples of \eqref{eq.gdef}, square roots $r=2$, degree $\deg p=3$

\begin{tabular}{ll}
$p(x)$ & \\
\hline
$1-4x^2-4x^3$ & A115962 \\
$1-2x-7x^2+8x^3$ & A098477 \\
$1-2x-3x^2+4x^3$ & A026569 \\
$1-2x-3x^2-4x^3$ & A191354 \\
$1-2x+x^2-4x^3$ & A098479 \\
$1-2x+x^2-8x^3$ & A098480 \\
$1-4x-8x^2-4x^3$ & A137635 \\
$1-4x+8x^3$ & A165431 \\
$1-4x+4x^3$ & A157004 \\
\end{tabular}

\end{exa}

\begin{exa}
Examples of \eqref{eq.gdef}, roots $r\neq 2$

\begin{tabular}{lll}
$p(x)$ & $r$ & \\
\hline
$1-4x$ & $2/3$ & A002457 \\
$1-8x$ & $2/3$ & A115902 \\
$1-36x$ & $6/11$ & A004998 \\
$1+9x+9x^3$ & $-3$ & A298308 \\
$1-9x-27x^3$ & $3$ & A095776 \\
\end{tabular}

\end{exa}

\section{Generalized Inverse Root}\label{sec.gisqr} 
The case with a non-trivial numerator polynomial $q(x)$ and denominator polynomial
$v(x)$ generalizes the content of Chapter \ref{sec.iroot}:
\begin{equation}
g(x) = \frac{q(x)}{v(x)\sqrt[r]{p(x)}},
\label{eq.qprootr}
\end{equation}
where expansion coefficients $q_n$ and $v_n$ are defined via
\begin{equation}
q(x) \equiv \sum_{n=0}^{\deg q} q_n x^n;\quad
v(x) \equiv \sum_{n=0}^{\deg v} v_n x^n.
\end{equation}
Multiply \eqref{eq.qprootr} by $v(x)$,
omitting the argument $x$ for brevity:
\begin{equation}
vg = qp^{-1/r}
\end{equation}
The derivative of this equation with respect to $x$ is
\begin{multline}
vg'+v'g = q'p^{-1/r}-\frac{1}{r}qp'p^{-1-1/r}\\
 = q' [qp^{-1/r}/v] \frac{v}{q}-\frac{1}{r}p'v \frac{qp^{-1-1/r}}{v} 
 = q' g\frac{v}{q}-\frac{1}{r}p'v \frac{1}{p}g.
\end{multline}
Multiplied by $rqp$ this is a first-order D-finite differential equation with polynomial
coefficients:
\begin{equation}
rqpv g' + (rqpv' - rpq'v +qp'v)g=0.
\label{eq.gfin}
\end{equation}
Assemble two auxiliary polynomials 
with coefficients $R_n$ 
and coefficients $Q_n$:
\begin{equation}
R(x) \equiv rqpv = \sum_{n=0}^{\deg R} R_n x^n;\quad
Q(x) \equiv rqpv' - rpq'v +qp'v = \sum_{n=0}^{\deg Q} Q_n x^n,
\label{eq.RQdef}
\end{equation}
such that \eqref{eq.gfin} reads
\begin{equation}
Rg' + Qg=0.
\label{eq.RQzero}
\end{equation}
\begin{rem}
\eqref{eq.qprootr} is a convolution of the sequence with generating function $q/v$
by the sequence with generating function \eqref{eq.gdef}. Compatible with Stoll's remark \cite{StollEJC18},
the product $R=rqpv$ indicates that the product of the degrees of the recurrences
of the convoluted sequences is a bound for the degree of the recurrence.
\end{rem}
\begin{rem}
If $p(x)$ is not a polynomial but a rational polynomial, the structure
is preserved if $R(x)$ and $Q(x)$ are multiplied by the denominators
that appear through the evaluation of $p$ and $p'$, so the format \eqref{eq.RQzero}
stays. The generic background is: replacing a generating function 
$g(x)\to \frac{q(x)}{v(x)}g(\frac{w(x)}{p(x)})$ with polynomials $q, v, w$ and $p$
preserves holonomicity (which is obvious because the chain rule applied to
the complicated generating function merely emits additional rational polynomial 
factors in comparison to the
differential equation of the simple $g(x)$). 
The binomial transformation \cite{BernsteinLAA226} of a P-finite sequence is 
again P-finite, for example.
\end{rem}
\begin{equation}
\sum_{i=0}^{\deg R} R_i x^i \sum_{j\ge 1} jg_jx^{j-1} + \sum_{i=0}^{\deg Q} Q_i x^i \sum_{j\ge 0}g_j x^j=0.
\end{equation}
\begin{multline}
\sum_{k=0}^{\deg R-1} \sum_{j=1}^{k+1} R_{k+1-j} jg_jx^k 
+\sum_{k\ge \deg R} \sum_{j=k+1-\deg R}^{k+1} R_{k+1-j} jg_jx^k  \\
+ \sum_{k=0}^{\deg Q-1} \sum_{j=0}^k Q_{k-j} g_j x^k
+ \sum_{k\ge \deg Q}  \sum_{j=k-\deg Q}^kQ_{k-j} g_j x^k
=0.
\end{multline}
Comparison of coefficients $[x^k]$ for sufficiently large $k$ on both sides
yields the P-finite recurrence
\begin{equation}
\sum_{j=k+1-\deg R}^{k+1} jR_{k+1-j} g_j
+ \sum_{j=k-\deg Q}^kQ_{k-j} g_j
=0,\quad k\ge \max(\deg R,\deg Q).
\end{equation}

Flip the direction in both $j$-sums:
\begin{equation}
\sum_{j=0}^{\deg R} (k+1-j)R_j g_{k+1-j}
+ \sum_{j=0}^{\deg Q} Q_j g_{k-j}
=0,\quad k\ge\max(\deg R,\deg Q).
\end{equation}
Substitute $k+1=n$:
\begin{equation}
\sum_{j=0}^{\deg R} (n-j)R_j g_{n-j}
+ \sum_{j=0}^{\deg Q} Q_j g_{n-j-1}
=0,\quad n>\max(\deg R,\deg Q).
\end{equation}
Replace $j\to j-1$ in the second term:
\begin{equation}
\sum_{j=0}^{\deg R} (n-j)R_j g_{n-j}
+ \sum_{j=1}^{\deg Q+1} Q_{j-1} g_{n-j}
=0,\quad n>\max(\deg R,\deg Q).
\end{equation}
In the case of \eqref{eq.RQdef} the degree of $Q$ is one less than the degree of $R$
because it contains one more derivative.
If we define coefficients $Q_n$ or $R_n$ to be zero if $n<0$ or $n$ larger than the degrees, this may be condensed as:
\begin{thm}\label{thm.P1st}
The P-finite recurrence of a sequence with the generating function \eqref{eq.qprootr}
is
\begin{equation}
\sum_{j=0}^{\deg R} [(n-j)R_j + Q_{j-1}]g_{n-j}
=0,\quad n>\deg R,
\label{eq.pfin1}
\end{equation}
where $R(x)$ and $Q(x)$ are the polynomials defined by the
sums and derivatives \eqref{eq.RQdef} of the three
polynomials $p(x)$, $q(x)$ and $v(x)$.
\end{thm}

\begin{rem}
To keep the recurrences simple, common polynomial factors of $x$ in
homogeneous differential equations like \eqref{eq.RQzero} should be eliminated (for example
by finding the greatest common divisor with the Euclidean algorithm)
before defining $R$ and $Q$.
\end{rem}

\begin{exa}
Examples of \eqref{eq.qprootr} with square roots $r=2$:

\begin{tabular}{llll}
$q(x)$ & $p(x)$ & $v$ & \\
\hline
$1-x$ & $1-6x+x^2$ & 1& A110170 \\
$1+x$ & $1-6x+x^2$ & 1& A241023 \\
$1-x$ & $1-6x+5x^2$ & 1& A085362 \\
$1-x$ & $1-2x-3x^2$ & 1& A025178 \\
$x(1+x)$ & $1-2x-3x^2$ & 1& A025565 \\
$1+2x$ & $1-4x^2$ & 1& A063886 \\
$1+x$ & $1+4x^2$ & 1& A128057 \\
$1$ & $1-4x$ & $1-x^2$ & A106188\\
$1$ & $1+4x$ & $1-4x$ & A091520\\
\end{tabular}

\end{exa}

We did not require that $r$ is an integer or positive. So formats like
$
g(x) = q(x)\sqrt{p(x)},\quad r=-2
$
or
$
g(x) = \frac{q(x)}{p^{3/2}(x)},\quad r=2/3
$
are also covered.

The format 
\begin{equation}
g=\sqrt{q(x)/p(x)}
\label{eq.gqp}
\end{equation}
 is reduced to the format \eqref{eq.qprootr} 
by multiplying numerator and denominator by $\sqrt{q(x)}$ such that the numerator is root-free.
The shortcut is:
\begin{thm}
The generating function \eqref{eq.gqp} obeys the differential equation
\begin{equation}
Rg'+Qg=0,
\label{eq.2term}
\end{equation}
with polynomials $R\equiv 2qp$ and $Q=qp'-q'p$, such that \eqref{eq.pfin1} applies.
\end{thm}

If an additive polynomial $w(x)$ appears on the right hand side like
\begin{equation}
g(x) = w(x) + \frac{q(x)}{v(x)\sqrt[r]{p(x)}},
\end{equation}
this modifies the coefficients $g_n$ for $n$ up to the degree
of the polynom $w(x)$. It delays the validity of \eqref{eq.pfin1}
to the point that all indices $j$ of the coefficients $g_j$ must
be larger than the degree of $w(x)$, so
$n$ in Theorem \ref{thm.P1st}
must be larger than the sum of $\tilde d$ and the polynomial degree of $w$.

\section{Generalized Inverse Root II} 
\subsection{Rooted Denominator}
The case with a non-trivial numerator polynomial $q(x)$ 
and denominator polynomials
$v(x)$ and
$w(x)$
generalizes the content further:
\begin{equation}
g(x) = \frac{q(x)}{w(x)+v(x)\sqrt[r]{p(x)}},
\label{eq.qprootr2}
\end{equation}
where expansion coefficients $q_n$, $v_n$ and $w_n$ are defined via
\begin{equation}
q(x) \equiv \sum_{n\ge 0} q_n x^n;\quad
v(x) \equiv \sum_{n\ge 0} v_n x^n;\quad
w(x) \equiv \sum_{n\ge 0} w_n x^n.
\end{equation}
With the strategy of Section \ref{sec.gisqr} one ends up with a differential
equation which contains terms proportional to $g^2$, which (apparently)
does not lead to recurrences with a finite number of terms.
If $r$ is a positive integer, $g$ is algebraic with $g^rv^rp=(q-gw)^r$.
So $g(x)$ is also holonomic \cite[Thm. 2.1]{StanleyEJC1}\cite{BanderierCPC24}.
Comtet's long division algorithm then yields the polynomial coefficients
of the associated linear differential equation of $g(x)$ \cite{ComtetLA10,CocklePM83,BostanISAAC07},
Final reshuffling of the coefficients of these polynomials 
by lowering the exponents
gives the P-finite recurrence.

Some progress can be made for square roots, $r=2$,
multiplying numerator and denominator of the fraction by $w-v\surd p$:
\begin{equation}
(w^2-v^2p)g = wq-vq\sqrt{p}.
\label{eq.tmp1}
\end{equation}
The first derivative is
\begin{multline}
(w^2-v^2p)g'+(w^2-v^2p)'g = (wq)'-(vq)'\sqrt{p}-vqp'\frac{1}{2\surd p}
\\
= (wq)'-\left[(vq)'+vq\frac{p'}{2p}\right]\sqrt{p}.
\end{multline}
Multiply by $2p$ to eliminate all denominators,
\begin{equation}
2p(w^2-v^2p)g'+2p(w^2-v^2p)'g
= 2p(wq)'-[2(vq)'p+vqp']\sqrt{p}.
\end{equation}
To eliminate the square root, multiply 
this equation by $vq$, multiply
\eqref{eq.tmp1} by $2(vq)'p+vqp'$, 
and subtract both equations
\begin{multline}
2pvq(v^2p-w^2)g'+2pvq(v^2p-w^2)'g 
+
[2(vq)'p+vqp'](w^2-v^2p)g 
\\
= -2pvq(wq)'
+ [2(vq)'p+vqp']wq.
\end{multline}

This is a first order differential equation with polynomial coefficients $R(x)$, $Q(x)$ and $H(x)$,
\begin{equation}
R(x)g'+Q(x)g=H(x),
\label{eq.RQH}
\end{equation}
where
\begin{equation}
R(x) \equiv 2pvq(v^2p-w^2) \equiv \sum_{j=0}^{\deg R} R_jx^j;
\label{eq.Rdef2}
\end{equation}
\begin{equation}
Q(x) \equiv -4pvqww'+2pq(pv^2+w^2)v'+vq(pv^2+w^2)p'-2pv(pv^2-w^2)q' \equiv \sum_{j=0}^{\deg Q} Q_jx^j;
\label{eq.Qdef2}
\end{equation}
\begin{equation}
H(x) \equiv -q^2(2pvw'-2wpv'-wvp').
\end{equation}
\begin{rem}
For $w=0$ the polynomials reduce to $H=0$, 
$R=2p^2v^3q$, 
$Q=pv^2(2pqv'+vqp'-2pvq')$.
The differential equation can be divided by the common factor $pv^2$ of $R$ and $Q$ because $H$ is zero,
and \eqref{eq.RQdef}--\eqref{eq.RQzero} emerge as a special case.
\end{rem}
\begin{rem}
Unlike \eqref{eq.RQH}, holonomic functions are defined to obey a differential
equation where no term such as $H(x)$ exists, which is not $g$ or a derivative of $g$.
That format is obtained by differentiating \eqref{eq.RQH} $d/dx$ $1+\deg H$ times,
such that $H$ disappears in the final higher-order differential equation \cite{StanleyEJC1}.
The philosophy in this paper is to keep the order of the differential
equation as low as possible, even if that means that the recurrence steps in 
belatedly due to the influence of the $H(x)$ on the low powers of $x$.
\end{rem}
The further reduction follows exactly the path of Section \ref{sec.gisqr},
paying attention to eliminate the early $a$-coefficients where the low, non-vanishing orders
of $H(x)$ may interfere:
\begin{thm}\label{thm.P2nd}
The P-finite recurrence of a sequence with the generating function \eqref{eq.qprootr2}
at $r=2$
is
\begin{equation}
\sum_{j=0}^{\deg R} [(n-j)R_j + Q_{j-1}]g_{n-j}
=0,\quad n-\deg R>\deg H,
\label{eq.P2nd}
\end{equation}
where $R(x)$ and $Q(x)$ are the polynomials defined by the
sums and derivatives \eqref{eq.Rdef2}--\eqref{eq.Qdef2} of the four
polynomials $p(x)$, $q(x)$, $v(x)$ and $w(x)$.
\end{thm}

\begin{exa}
Examples of \eqref{eq.qprootr2}:

\begin{tabular}{llllll}
$p$ & $q$ & $v$ & $w$ & $r$ & \\
\hline
$1-2x-3x^2$ & 1 & $1+x$ & $-x$ & 2 & A116394 \\
\end{tabular}
\end{exa}

\subsection{Rooted Numerator}\label{sec.galge}
If the generating function has the shape
\begin{equation}
g(x)=\frac{w(x)+v(x)\sqrt[r]{p(x)}}{q(x)}
\label{eq.qprootr3}
\end{equation}
multiply this equation by $q$, derive it, and eliminate $\sqrt[r]{p}$
with the aid of the previous equation:
\begin{equation}
qg'+q'g=w'+v'\sqrt[r]{p}+\frac{1}{r} v\frac{p'}{p}\sqrt[r]{p}
=w'+\left[v'+\frac{vp'}{rp}\right]\frac{qg-w}{v};
\end{equation}
Multiply by $rpv$ to obtain an equation which fits \eqref{eq.RQH}, this time with
\begin{equation}
R(x) \equiv rpqv = \sum_{n=0}^{\deg R}R_n x^n;
\label{eq.Rdef3}
\end{equation}
\begin{equation}
Q(x) \equiv rp(q'v-qv')-vp'q = \sum_{n=0}^{\deg Q} Q_n x^n;
\label{eq.Qdef3}
\end{equation}
\begin{equation}
H(x) \equiv rp(w'v-wv')-vp'w;
\label{eq.Hdef3}
\end{equation}
\begin{thm}
The P-finite recurrence of a sequence with the generating function \eqref{eq.qprootr3}
is
given by \eqref{eq.P2nd}
where $R(x)$ and $Q(x)$ are the polynomials defined by the
sums and derivatives \eqref{eq.Rdef3}--\eqref{eq.Qdef3} of the four
polynomials $p(x)$, $q(x)$, $v(x)$ and $w(x)$.
\end{thm}
This formula and Section \ref{sec.iroot} cover Callan's generating functions \cite{CallanJIS10}.

The generating functions of the form
\begin{equation}
g=\frac{u(x)}{w(x)}+\frac{q(x)}{v(x)\sqrt[r]{p(x)}}
\end{equation}
with polynomials $p(x)$, $v(x)$, $u(x)$ and $w(x)$
are also covered by the form \eqref{eq.qprootr3}
because they can be rewritten as
\begin{equation}
g=\frac{u(x)v(x)p(x)+w(x)q(x)p^{1-1/r}(x)} {w(x)v(x)p(x)}.
\end{equation}
This allows to find P-recurrences of sequence which are sums
of C-finite sequences represented by $u(x)/w(x)$ and sequences represented by $q(x)/[v(x)\sqrt[r]{p(x)}]$,
such as transiting from \cite[A026375]{sloane} to \cite[A242586]{sloane}.

The form with a common square root in numerator and denominator
is also in this class:
\begin{multline}
g=\frac{u(x)\sqrt{p(x)}+w(x)}{v(x)\sqrt{p(x)}+q(x)}
\\
=
\frac{u(x)v(x)p(x)-q(x)w(x)+[w(x)v(x)-q(x)u(x)]\sqrt{p(x)}}{v^2(x)p(x)-q^2(x)}
.
\end{multline}

\subsection{Orthogonal Polynomials}
Some orthogonal polynomials have generating functions which are in our category
of rational square root
expressions \cite[\S 22.9]{AS}. If the argument $x$ of these orthogonal polynomials
is kept fixed and their degree $n$ used as the index of the $g(n)$,
their well-known 3-term recurrences appear \cite[22.7]{AS}.

\section{Generalized Hypergeometric Function} 
\subsection{Power Series}
The Generalized Hypergeometric Functions $_pF_q(x)$ 
with a set of constant ``numerators'' $\{\alpha\}_p$
and ``denominators'' $\{\beta\}_q$ are another special case
with simple P-finite recurrences \cite{SlaterHyp}:
\begin{equation}
g(x) = x^t{} _pF_q(\{\alpha\}_p; \{\beta\}_q; x^r/c)
=\sum_{n\ge 0} \frac{\prod_{i=1}^p (\alpha_i)_n}{\prod_{j=1}^q (\beta_j)_n}\frac{x^{t+rn}}{n!c^n},
\label{eq.pFq}
\end{equation}
where $(\cdot)_{\cdot .}$ are Pochhammer symbols \cite[(13.1.2)]{AS}
\begin{equation}
(\alpha)_n\equiv \frac{\Gamma(\alpha+n)}{\Gamma(\alpha)} =\alpha(\alpha+1)(\alpha+2)\cdots (\alpha+n-1).
\end{equation}
The non-vanishing coefficients of the power series are
\begin{equation}
g_{rn+t} =
\frac{\prod_{i} (\alpha_i)_n}{\prod_j (\beta_j)_n}\frac{1}{n!c^n}.
\end{equation}
The associated P-finite 2-term recurrence is
\begin{equation}
c(n+1) \prod_j(\beta_j+n) g_{rn+r+t}
=
\prod_{i} (\alpha_i+n) g_{rn+t}.
\end{equation}
\begin{exa}
The generating function
\begin{equation}
(1-x)^\alpha = {} _1F_0(-\alpha;;x)
\end{equation}
is a borderline case between \eqref{eq.gdef} and \eqref{eq.pFq}. Also
\begin{eqnarray}
\arcsin x &=& x{}_2F_1(\frac12,\frac12; \frac32; x^2),\\
J_\nu(z)&=&\frac{(z/2)^\nu}{\Gamma(\nu+1)}{}_0F_1(\nu+1;-z^2/4),
\end{eqnarray}
and
\begin{equation}
\ln(1+x) = x{}_2F_1(1,1;2;-x)
\end{equation}
fit in here \cite[\S 0.7.2]{Bronstein2}\cite{KoepfJSC13}.
\end{exa}

\subsection{Sequence Index in Parameters}
There is another family of recurrences associated with hypergeometric functions.
If the sequence entries $g(n)$ are hypergeometric functions of some constant
argument $x$ where the index $n$ appears in the parameters such that  $\alpha_p$
or $\beta_p$ are polynomials of $n$, the 3-term contiguous relations 
of the Gaussian Hypergeometric Functions \cite[\S 15.2]{AS},
the 3-term recurrences of the Confluent Hypergeometric Functions \cite[\S 13.4]{AS}\cite{RakhaCMA61}
and the contiguous relations for other pairs of $(p,q)$ \cite{RainvilleBAMS51,BaileyPGMA2,VidunasJCAM153}
provide the P-finite recurrences for $g(n)$.

A side aspect is that for the \emph{terminating} hypergeometric functions in that case, i.e., polynomials
of $n$, the P-finite recurrences are even C-finite.

\section{Nested Roots} \label{sec.nroot} 
\subsection{Elliptic Integrals}
Let
\begin{equation}
g(x)=\sqrt[r]{f(x)}
\end{equation}
be a generating function with a discriminant function $f(x)$. Then
\begin{equation}
g' = \frac{1}{r}\frac{f'}{f}g;
\label{eq.nest1}
\end{equation}
Suppose also that $f(x)$ obeys a D-finite first-order 
differential equation of the form
\begin{equation}
L(x)f'+H(x)f=0,
\end{equation}
with polynomials $L(x)$ and $H(x)$, which implies that $\log f(x)=-\int \frac{H(x)}{L(x)}dx$ are Elliptic Integrals. 
Then  \eqref{eq.nest1} becomes
\begin{equation}
rL(x)g' + H(x)g =0.
\end{equation}

This means if we have obtained a D-finite first-order differential equation for a generating
function, the D-finite differential equation equation for a generating function
that is some (fractional) power of the original generating function is inherited
(apart from the factor $r$) by the derived generating function \cite{KoepfJSC13}.
For squares of hypergeometric functions see e.g. Chaundy
and Vidunas \cite{ChaundyQJM14,VidunasRam26}.

\subsection{Polynomial Discriminants}
Let
\begin{equation}
g(x)=\sqrt[r]{w(x)+\sqrt{p(x)}}
\label{eq.gnroot}
\end{equation}
be a generating function with polynomials $w(x)$ and $p(x)$.

The first and second derivatives of \eqref{eq.gnroot} are
\begin{equation}
g'
=\frac1r (w'+\frac12 \frac{p'}{\sqrt p})(w+\sqrt{p})^{1/r-1};
\end{equation}
\begin{equation}
g''
=\frac1r\left[(\frac1r-1) (w'+\frac12 \frac{p'}{\sqrt p})^2(w+\sqrt p)^{\frac1r -2}
+(w''+\frac12 \frac{p''}{\sqrt p}-\frac14 \frac{{p'}^2}{p^{3/2}})(w+\sqrt p)^{\frac1r -1}
\right].
\end{equation}

With the ansatz
\begin{equation}
T(x) g'' + R(x) g' + Q(x) g=0
\label{eq.TRQdeq}
\end{equation}
we \emph{assume} that this generating function is D-finite with three
polynomials $T(x)$, $R(x)$ and $Q(x)$.  This requires
\begin{multline}
T(x)
\frac1r\left[(\frac1r-1) (w'+\frac12 \frac{p'}{\sqrt p})^2(w+\sqrt p)^{\frac1r -2}
+(w''+\frac12 \frac{p''}{\sqrt p}-\frac14 \frac{{p'}^2}{p^{3/2}})(w+\sqrt p)^{\frac1r -1}
\right]
\\
+ R(x)
\frac1r (w'+\frac12 \frac{p'}{\sqrt p})(w+\sqrt{p})^{1/r-1}
+ Q(x)
(w+\sqrt{p})^{1/r}
=0,
\end{multline}
and multiplied by $r^2(w+\surd p)^{2-\frac1r}$
\begin{multline}
T(x)
\left[(1-r) (w'+\frac12 \frac{p'}{\sqrt p})^2
+r(w''+\frac12 \frac{p''}{\sqrt p}-\frac14 \frac{{p'}^2}{p^{3/2}})(w+\sqrt p)
\right]
\\
+ R(x)
r (w'+\frac12 \frac{p'}{\sqrt p})(w+\sqrt{p})
+ Q(x)
r^2(w+\sqrt{p})^2
=0.
\end{multline}
Expanding squares and products this reads
\begin{multline}
T(x)
[
(1-r){w'}^2+(1-r)w'\frac{p'}{\sqrt p}+(\frac14-\frac{r}{2}) \frac{{p'}^2}{p}
\\
+rw''w+\frac12 r w\frac{p''}{\sqrt p}-\frac14 rw\frac{{p'}^2}{p^{3/2}}
+rw''\sqrt p +\frac12 rp''
]
\\
+ R(x)
(rw'w+\frac12 rw \frac{p'}{\sqrt p}
+rw'\sqrt p +\frac12 rp')
+ Q(x)
(r^2w^2+2r^2w\sqrt{p}+r^2p)
=0.
\end{multline}

This has the structure
\begin{equation}
T(x)[\alpha_1(x)+\frac{1}{p^{3/2}}\alpha_4(x)]
+R(x)[\alpha_2(x)+\frac{1}{p^{3/2}}\alpha_5(x)]
+Q(x)[\alpha_3(x)+\frac{1}{p^{3/2}}\alpha_6(x)]=0,
\label{eq.TRQzero}
\end{equation}
where 6 $\alpha$-coefficients are rational polynomials in $x$ defined as
\begin{eqnarray}
\alpha_1(x) &\equiv& 
(1-r){w'}^2+(\frac14-\frac{r}{2}) \frac{{p'}^2}{p} +rw''w+\frac12 rp''
;
\\
\alpha_2(x) &\equiv& 
rw'w +\frac12 rp'
;
\\
\alpha_3(x) &\equiv& 
r^2w^2+r^2p
;
\\
\alpha_4(x) &\equiv& 
(1-r)w'p'p
+\frac12 r wp''p-\frac14 rw{p'}^2
+rw''p^2
;
\\
\alpha_5(x) &\equiv& 
\frac12 rw p'p
+rw'p^2 
;
\\
\alpha_6(x) &\equiv& 
2r^2wp^2
.
\end{eqnarray}
Instead of solving \eqref{eq.TRQzero} in general we continue with a separation
ansatz, where 
the components which depend on $1/p^{3/2}$ and do not depend on it are individually
zero:
\begin{eqnarray}
T(x)\alpha_1(x)
+R(x)\alpha_2(x)
+Q(x)\alpha_3(x)&=&0 ;\\
\wedge\quad
T(x)\alpha_4(x)
+R(x)\alpha_5(x)
+Q(x)\alpha_6(x)&=&0.
\end{eqnarray}
In the language of 3-dimensional vector algebra, the vector $(T,R,Q)$
is orthogonal to the vector $(\alpha_1,\alpha_2,\alpha_3)$ as well
as orthogonal to the vector $(\alpha_4,\alpha_5,\alpha_6)$, so it is the vector cross product
of the two $\alpha$-vectors:
\begin{eqnarray}
T(x)&=& \alpha_2\alpha_6-\alpha_3\alpha_5;
\label{eq.Talpha}
\\
R(x)&=& \alpha_3\alpha_4-\alpha_1\alpha_6;
\\
Q(x)&=& \alpha_1\alpha_5-\alpha_2\alpha_4.
\label{eq.Qalpha}
\end{eqnarray}
\begin{rem}
This construction of a $d$-dimensional vector $V$ that is orthogonal to 
$d-1$ vectors $v^{(1)}, v^{(2)},\ldots, v^{(d-1)}$
such that $\sum_{i=1}^d V_i v^{(1)}_i = \sum_{i=1}^d V_i v^{(2)}_i =\cdots =0$
carries over to more than 3 dimensions \cite{ShawIJMEST18}: the $i$-th component of
$V$ is the tensor sum-product (determinantal mix)
of the $\epsilon$-Tensor (parity
of the permutation of its indices) with the product of the components of the $v$-vectors:
\begin{equation}
V_i = \sum_{j=1}^d \sum_{k=1}^d\cdots \sum_{m=1}^d \epsilon_{ijk\cdots m} v^{(1)}_j v^{(2)}_k \cdots v^{(d-1)}_m.
\end{equation}
\label{rem.alg}
\end{rem}
Insertion of the $\alpha$-terms
and multiplying $T$, $R$ and $Q$ with a common
factor $8/r$ yields:
\begin{equation}
T= 
4r^2(-w^2+p)p(-2w'p+wp')
;
\label{eq.Tnroot}
\end{equation}
\begin{multline}
R= 
-2r(-4w^2w'p'p+4rw^2w'p'p-2rw^3p''p+rw^3{p'}^2+4rw^2w''p^2-4p^2w'p'+4rp^2w'p'\\
+2rp^2wp''
+5rwp{p'}^2-4rp^3w''+8wp^2{w'}^2-8rwp^2{w'}^2-2wp{p'}^2)
;
\end{multline}
\begin{multline}
Q=
-4w'^2wp'p+8{w'}^3p^2+4{w'}^2rwp'p-8{w'}^3rp^2-{p'}^3w-6{p'}^2pw'
+3{p'}^3rw
\\
+8{p'}^2prw'+4rw''w^2p'p+4rp''w'
p^2-4rw'w^2p''p+2rw'w^2{p'}^2-4rp'w''p^2
.
\label{eq.Qnroot}
\end{multline}
Working backwards through the logic shows that the ansatz \eqref{eq.TRQdeq} is
indeed satisfied.

\begin{thm}
The coefficients of the generating function \eqref{eq.gnroot} obey the P-finite recurrence
\begin{equation}
\sum_{j\ge 0} \left[(n-j)(n-j-1)T_j+(n-j)R_{j-1}+Q_{j-2}\right]g_{n-j} 
=0,
\label{eq.TRQrec}
\end{equation}
where $T=\sum_{n\ge 0}T_nx^n$, $R=\sum_{n\ge 0} R_nx^n$
and $Q=\sum_{n\ge 0}Q_nx^n$ are the polynomials \eqref{eq.Tnroot}--\eqref{eq.Qnroot}.
\end{thm}

If the polynomials of $n$ in front of the $g_{n-j}$ are given
in the standard basis of powers of $n$, the coefficients $T_j$,
$R_j$ and $Q_j$ are easily recovered by accumulating \cite[24.1.4]{AS}\cite{SalvyTOMS20}
\begin{equation}
n^t = \sum_{m=0}^t {\mathcal S}_t^{(m)}n(n-1)(n-2)\cdots (n-m+1),
\label{eq.ntoStir}
\end{equation}
where $\mathcal S$ are the Stirling Numbers of the Second Kind.

\subsection{Degenerate cases}
$T(x)$ of \eqref{eq.Tnroot} is zero if $w^2=p$ or $wp'=2w'p$, and the simpler \eqref{eq.RQzero} applies;
the P-recurrence only involves first-degree polynomials.
The case $w^2=p$ is not interesting: then $g=\sqrt[r]{2w}$ has the format \eqref{eq.gdef} and
would be treated accordingly.
The case $wp'=2w'p$ means $2w'/w=p'/p$, therefore $2\ln w=\ln p+C$, 
therefore  $\ln w^2=\ln p+C$, therefore $w^2=Cp$, and again \eqref{eq.gdef} is the
underlying format.

\subsection{Deeper Nests}
The roots may be nested deeper where
\begin{equation}
g(x)=\sqrt[r_k]{w_k(x)+\sqrt[r_{k-1}]{w_{k-1}(x)+\cdots \sqrt[r_1]{w_1(x)}}},
\end{equation}
were the $w_i(x)$ are polynomials and where the $r_i$ are integers. 
[A first approach to obtain the coefficients $g_n$ numerically is to regard
this as the composition of roots \cite{BrentJACM25}.]

$g(x)$ is an algebraic
function: take the $r_k$-th power of both sides, move $w_k$ to the left side,
take the $r_{k-1}$-st power, move $w_{k-1}(x)$ to the left side and so on
to establish its algebraic equation.  The highest power is $g^{r_kr_{k-1}\cdots r_1}(x)$.
The generic strategy of Section \ref{sec.galge} establishes a recurrence.

\section{Exponentials} 
\subsection{Exponential of Root}

The class of generating functions
\begin{equation}
g(x)=\exp[w(x)\pm\sqrt{p(x)}]
\label{eq.expsqrt}
\end{equation}
with polynomials $w(x)$ and $p(x)$
has similar regenerative properties as the nested roots:
\begin{equation}
g'=\left(w'\pm \frac12 p'p^{-1/2}\right)\exp[w\pm \sqrt p];
\end{equation}
\begin{equation}
g''=
\left(w''\pm\frac12 p''p^{-1/2}\mp\frac14 {p'}^2p^{-3/2}\right)\exp[w\pm\sqrt p]
+\left(w'\pm\frac12 p'p^{-1/2}\right)^2\exp[w\pm\sqrt p].
\end{equation}
\begin{rem}
If $p(x)=0$, the subsequent sections are superfluous because then $g'-w'g=0$,
which is a special case of \eqref{eq.RQzero}, so the recurrence \eqref{eq.pfin1} applies.
\end{rem}
The same procedure as in Section \ref{sec.nroot} unfolds:
\begin{multline}
T
[\left(w''\pm\frac12 p''p^{-1/2}\mp\frac14 {p'}^2p^{-3/2}\right)\exp[w\pm\sqrt p]
+\left(w'\pm\frac12 p'p^{-1/2}\right)^2\exp[w\pm\sqrt p]]\\
+R
\left(w'\pm\frac12 p'p^{-1/2}\right)\exp[w\pm\sqrt p]
+Q
\exp[w\pm\sqrt p]=0.
\end{multline}
\begin{multline}
T
[w''\pm\frac12 p''p^{-1/2}\mp\frac14 {p'}^2p^{-3/2}
+\left(w'\pm\frac12 p'p^{-1/2}\right)^2]
+R
\left(w'\pm\frac12 p'p^{-1/2}\right)
+Q= 0.
\end{multline}
\begin{multline}
T
[w''\pm\frac12 p''p^{-1/2}\mp\frac14 {p'}^2p^{-3/2}
+{w'}^2\pm w'p'p^{-1/2}+\frac14 {p'}^2\frac{1}{p} ]
+R
\left(w'\pm\frac12 p'p^{-1/2}\right)
+Q= 0.
\end{multline}
\begin{rem}
If $\surd p(x)$ is replaced by $\sqrt[r]p$ in the generating function \eqref{eq.expsqrt},
terms proportional to $p^{1/r}$, proportional to $p^{2/r}$ and not
related to $p^{1/r}$ appear, and a higher-dimensional ansatz following
Remark \ref{rem.alg} is needed to disentangle the three algebraic branches.
\end{rem}
The coefficients matching \eqref{eq.TRQzero} are
\begin{eqnarray}
\alpha_1 &=& w''+{w'}^2 +\frac14 \frac{{p'}^2}{p};
\\
\alpha_2 &=& w';
\\
\alpha_3 &=& 1;
\\
\alpha_4 &=& \pm\frac12 p'' p\mp\frac14 {p'}^2 \pm w'p'p;
\\
\alpha_5 &=& \pm\frac12 p'p;
\\
\alpha_6 &=& 0.
\label{eq.alpha6Qexpsqrt}
\end{eqnarray}
With these we gather three polynomials \eqref{eq.Talpha}--\eqref{eq.Qalpha} 
and multiply $T$, $R$ and $Q$ by the common factor 8 to maintain integer coefficients:
\begin{eqnarray}
T(x) &=&  \mp 4p'p =\sum_{n\ge 0}T_nx^n;
\label{eq.Texpsqrt}
\\
R(x) &=& \pm 2(2p''p-{p'}^2+4w'p'p) =\sum_{n\ge 0}R_nx^n;\\
Q(x) &=& \pm(4p'pw''-4p'p{w'}^2+{p'}^3-4w'p''p+2w'{p'}^2)=\sum_{n\ge 0}Q_nx^n.
\label{eq.Qexpsqrt}
\end{eqnarray}

\begin{thm}\label{thm.Qexpsqrt}
The coefficients of the generating function \eqref{eq.expsqrt} obey the P-finite recurrence
\eqref{eq.TRQrec}
where $T(x)$, $R(x)$
and $Q(x)$ are the polynomials \eqref{eq.Texpsqrt}--\eqref{eq.Qexpsqrt}.
\end{thm}

\subsection{Exponential Times Arithmetic} 

If the generating function is of the kind
\begin{equation}
g(x) = \exp[q(x)/v(x)]\frac{1}{\sqrt[r]{p(x)}}
\label{eq.ep}
\end{equation}
with polynomials $p(x)$, $q(x)$ and $v(x)$,
then by the product and chain rules
\begin{equation}
g' = \left(\frac{q'}{v}-\frac{qv'}{v^2}\right)\exp(q/v)\frac{1}{\sqrt[r]{p}}
-\frac1r \exp(q/v)\frac{p'}{p^{1+1/r}}
= (\frac{q'}{v}-\frac{qv'}{v^2}) g -\frac1r \frac{p'}{p}g .
\end{equation}
Multiplication with $rpv^2$ yields the first order differential equation
\begin{equation}
rpv^2g' = [rp(q'v-qv')-v^2p']g.
\end{equation}
So we face 
\eqref{eq.RQzero} 
but this time with the polynomials
\begin{equation}
R(x)\equiv rpv^2,\quad Q(x)\equiv v^2p'-rp(q'v-qv').
\end{equation}
and apply the recurrence 
\eqref{eq.pfin1}.

\section{Logarithm of Rational} 
If the generating function is of the kind
\begin{equation}
g(x)= \log[q(x)/v(x)]
\end{equation}
with polynomials $p(x)$ and $q(x)$, then
\begin{equation}
g'(x)= \frac{v}{q}[\frac{q'}{v}-\frac{qv'}{v^2}].
\end{equation}
This is a special case of \eqref{eq.RQH} substituting
\begin{equation}
R(x)= qv;\quad Q(x)=0;\quad H(x)= vq'-qv'
\end{equation}
and obeys the recurrence \eqref{eq.P2nd}. If in addition the Wronskian $H$ is zero
(i.e., if the two polynomials have a common factor),
this simplifies furthermore to \eqref{eq.RQzero}.

\begin{rem}
The logarithm of the ratio is the difference between the logarithms of numerator and denominator.
These logarithms are separately holonomic, and by the closure property one may also generate
the P-finite recurrence from the recurrences of the two terms. The technique to combine
the two recurrences for a Hadamard sum has been demonstrated by Mallinger \cite[Thm. 1.4.3]{MallingerMsc}
and Kauers \cite{KauersCAQTF}.
\end{rem}

\section{Summary}
We validated a set of P-recurrences
of sequences which involve generating function with
roots. 

\appendix
\section{Inhomogeneous P-finite}
If the sequence $a$ obeys a P-finite recurrence with a polynomial $I(n)$,
\begin{equation}
\sum_{j\ge 0} P_j(n) g_{n-j} + I(n)=0,
\label{eq.In}
\end{equation}
it can be rewritten as a homogeneous P-finite recurrence by shifting the index by 1:
\begin{equation}
\sum_{j\ge 0} P_j(n-1) g_{n-j-1} + I(n-1)=0,
\label{eq.In1}
\end{equation}
multiplying \eqref{eq.In} by $I(n-1)$ and multiplying \eqref{eq.In1} by $I(n)$
and subtracting both equations.
This results in a recurrence which is one longer than \eqref{eq.In} and
has polynomial coefficients of degrees which are the sum of the degrees in \eqref{eq.In}
and the degree of $I(n)$.

\section{Exponential Generating Functions}
A generating function $g$ may be interpreted as an \emph{ordinary} generating
function for a sequence $g(n)$ and at the same time as an \emph{exponential} generating
function for a sequence $b(n)$:
\begin{equation}
g(x) = \sum_{n\ge 0} g_n x^n = \sum_{n \ge 0} b_n \frac{x^n}{n!}.
\end{equation}
If the $g_n$ obey a P-finite recurrence with polynomials $P$ and $J+1$ terms,
\begin{equation}
\sum_{j= 0}^J P_j(n) g_{n-j} =0,
\end{equation}
substituting $g_n = b_n/n!$ and multiplication of the recurrence with $n!$ yields
an associated P-finite recurrence for the $b_n$ \cite{StanleyEJC1}:
\begin{equation}
\sum_{j=0}^J (n-j+1)_j P_j(n) b_{n-j} =0.
\end{equation}
So the polynomials in the P-finite recurrence of the $b$-terms of the \emph{exponential} generating
function are the polynomials of the $a$-terms of the \emph{ordinary} generating  function
multiplied by first-order polynomials with can be represented as Pochhammer symbols.
\begin{rem}
An equivalent match applies to logarithmic generating functions $g(x)=\sum_{n\ge 0} c_n x^n/n$.
\end{rem}

\begin{exa}
If we reinterpret \eqref{eq.qprootr} as an exponential generating function and multiply the
coefficients of \eqref{eq.pfin1} with the Pochhammer symbols, some OEIS sequences are covered:

\begin{tabular}{lllll}
$p(x)$ & $q(x)$ & $v(x)$ & $r$ &\\
\hline
$1-4x+x^2$ & 1& 1& 2& A285199 \\
$1-8x+x^2$ & 1& 1& 2& A006438 \\
$1+2x+4x^2$ & 1& 1& 2& A182827 \\
$1-2x-2x^2$ & 1& 1& 2& A098460 \\
$1-2x-3x^2$ & 1& 1& 2& A098461 \\
$1-10x$ & 1& 1& 10& A144773 \\
\end{tabular}
\end{exa}

\begin{exa}
If we reinterpret \eqref{eq.ep} as an exponential generating function
and multiply the coefficients of \eqref{eq.pfin1} with the Pochhammer symbols, additional OEIS sequences are covered:

\begin{tabular}{lllll}
$p(x)$ & $q(x)$ & $v(x)$ & $r$ \\
\hline
$1$ & $x$ & $1-x$ & $1$ & A000262 \\
$1-2x$ & $x$ & 1 & $-2$ & A055142 \\
$1-4x$ & $x$ & 1 & $2$ & A052143 \\
$1-8x$ & $8x$ & 1 & $8$ & A094935 \\
$1-7x$ & $7x$ & 1 & $7$ & A094911 \\
$1$ & $x(1+x)$ & $1-x-x^2$ & $1$ & A345075 \\
$1$ & $3x(2+x)/2$ & $1$ & $1$ & A335819 \\
$1-x$ & $-x(2+x)$ & $1$ & $1/2$ & A335595 \\
$1$ & $3x(2+x)/2$ & $1$ & $1$ & A335819 \\
$1+x$ & $x$ & $1-x$ & $1$ & A331725 \\
$1$ & $x$ & $(1+x)^2$ & $1$ & A318215 \\
$6+x(3+x)(6+x)$ & $-x$ & $1$ & $1$ & A302908 \\
$1-4x$ & $-2x$ & $1$ & $1$ & A296660 \\
$1$ & $x+x^2-x^3/6$ & $1$ & $1$ & A200380\\
\end{tabular}

\end{exa}

\section{Reduction Of the Number of Terms} 
In \cite[A122877]{sloane} the generating function
\begin{equation}
g = \frac{1-2x-3x^2-(1-x)\sqrt{1-2x-7x^2}}{8x^3}
\end{equation}
matches \eqref{eq.qprootr3} with polynomials $q=8x^3$, $w=1-2x-3x^2$,
$v=-(1-x)$, $p=1-2x-7x^2$ and $r=2$, such that
the differential equation \eqref{eq.RQH} is set up with
$R=-16x^3(1-x)(1-2x-7x^2)$, $Q=-16x^2(3-7x-11x^2+7x^3)$
and $H=-64x^3$ defined in \eqref{eq.Rdef3}--\eqref{eq.Hdef3}. 
The greatest common factor $-16x^2$ of $R$, $Q$ and $H$ can be dropped
in the differential equation:
\begin{equation}
x(1-x)(1-2x-7x^2)g' +(3-7x-11x^2+7x^3) g= 4x.
\label{eq.appdg1}
\end{equation}
Only the indices $j=1$--$4$ contribute to the
recurrence \eqref{eq.P2nd}, so the generating function supports a 4-term recurrence:
\begin{equation}
(n+3)g_n-(3n+4)g_{n-1}-(5n+1)g_{n-2}+7(n-2)g_{n-3}=0
\label{eq.exa4term}
\end{equation}
with first degree polynomials.
Differentiating of 
\eqref{eq.RQH} 
yields a second order 
differential equation
\begin{equation}
Rg''+(R'+Q)g'+Q'g=H',
\label{eq.RQzero2}
\end{equation}
here
\begin{equation}
x(1-x)(1-2x-7x^2)g'' +(1-x)(4-9x-35x^2)g'+(-7-22x+21x^2) g= 4.
\end{equation}
The number of the terms in the recurrence derived from the first-order differential
equation is based on:
\begin{itemize}
\item
The factor $R$ contributes powers $x^0$\ldots $x^{\deg R}$;
$g'$ represents $\sum ng_nx^{n-1}$, so the product has powers $x^{n-1}$
up to $x^{n-1+\deg R}$
\item
The factor $Q$ contributes powers $x^0$\ldots $x^{\deg Q}$;
$g$ represents $\sum g_nx^n$, so the product has powers $x^n$
up to $x^{n+\deg Q}$.
\end{itemize}
The range of powers is $x^{n-1}$ up to the larger of $x^{n-1+\deg R}$ or $x^{n+\deg Q}$,
and the spread of exponents determines the number of coupled $a$-coefficients.
The equivalent analysis of the second-order differential equation (using 
$g''\equiv \sum n(n-1)g_nx^{n-2}$) shows that the exponents have
been decremented by one, but the spread of exponents remains the same.
[This preservation remains valid, even if some lower coefficients $R_n$ vanish,
like in our example where $R_0=0$.]
The numbers of terms in the P-recurrences derived from 
\eqref{eq.RQH} and \eqref{eq.RQzero2} are the same. 

The \emph{penalty}
in \eqref{eq.RQzero2}, induced by $g''\sim \sum n(n-1)g_n$, is that the
polynomials in the P-recurrences are of degree 2, not 1.
\emph{However}, if $\sum_n R_n=0$ [equivalent: a factor $1-x$ in 
the factorization of $R(x)$],
the contribution of the $n^2$ terms in the recurrence derived from \eqref{eq.RQzero2}
vanishes. In that circumstance the P-recurrence from \eqref{eq.RQzero2} also has
coefficients which are polynomials of \emph{first} degree. 
In the synoptical view on both recurrences
of the same number of terms and the same polynomial degrees, one may
multiply each recurrence with the polynomial in front of $g_n$ of the other recurrence,
subtract both, to obtain a recurrence with one term less and with
polynomial coefficients
with a degree which is the sum of the individual degrees.

In the example considered here, the requirement on $R(x)$ is fulfilled,
and besides \eqref{eq.exa4term} there is a 3-term recurrence
\begin{equation}
-(n+3)(n-1)g_n+n(2n+1)g_{n-1}+7n(n-1)g_{n-2}=0
\end{equation}
with quadratic polynomials.

\begin{rem}
The derivative of D-finite differential equations with polynomial coefficients
yields differential equations of higher order, equivalent to P-recurrences with
polynomials of higher degrees, and potentially of smaller length. 
We take the stand that recurrences derived
from differential equations of lower order are preferable, even if the number of
terms in the P-recurrences (the \textit{length} of the recurrences) is larger,
because the step from the P-recurrences to the D-equation plus differentiation
is straight forward, whereas the opposite direction (one integration of the D-equation)
may be difficult and introduces further constants.
\end{rem}

\bibliographystyle{amsplain}
\bibliography{all}

\end{document}